\documentclass[journal]{IEEEtran}
\ifCLASSINFOpdf
\else
\fi
\usepackage[latin9]{inputenc}
\usepackage{color}
\usepackage{array}
\usepackage{graphicx}

\usepackage{tocloft}
\usepackage[caption=false,font=footnotesize]{subfig}

\cftpagenumbersoff{figure}

\begin{document}
\title{\textcolor{black}{ Artificial Intelligence for Dynamical Systems in Wireless Communications: Modeling for the Future}}
\author{Harun Siljak,~\IEEEmembership{Senior Member,~IEEE,} Irene Macaluso,
and~Nicola Marchetti,~\IEEEmembership{Senior Member,~IEEE}\thanks{Manuscript received ../../.... The authors are with the CONNECT Centre,
Trinity College Dublin, Ireland (emails: \{harun.siljak, macalusi,
nicola.marchetti\}@tcd.ie). This work was supported in part by a research
grant from Science Foundation Ireland (SFI), co-funded under the European
Regional Development Fund under Grant Numbers 13/RC/2077 and 13/RC/2077\_P2, as well as
the European Union\textquoteright s Horizon 2020 research and innovation
programme under the Marie Sk\l odowska-Curie grant agreement No 713567.}}
\maketitle
\begin{abstract}
\textcolor{black}{Dynamical systems are no strangers in wireless communications.
Our story will necessarily involve chaos, but not in the terms secure
chaotic communications have introduced it: we will look for the chaos,
complexity and dynamics that already exist in everyday wireless communications.
We present a short overview of dynamical systems and chaos before
focusing on the applications of dynamical systems theory to wireless
communications in the past 30 years, ranging from the modeling on
the physical layer to different kinds of self-similar traffic encountered
all the way up to the network layer. The examples of past research
and its implications are grouped and mapped onto the media layers
of ISO OSI model to show just how ubiquitous dynamical systems theory
can be and to trace the paths that may be taken now. When considering
the future paths, we argue that the time has come for us to revive
the interest in dynamical systems for wireless communications. It
did not happen already because of the big question: can we afford
observing systems of our interest as dynamical systems and what are
the trade-offs? The answers to these questions are dynamical systems
of its own: they change not only with the modeling context, but also
with time. In the current moment the available resources allow such
approach and the current demands ask for it. Reservoir computing,
the major player in dynamical systems-related learning originated
in wireless communications, and to wireless communications it should
return.}
\end{abstract}


\section{\textcolor{black}{Introduction}}

\textcolor{black}{Everything is a dynamical system, depending on how
you define everything (and how you define a dynamical system). In
the realm of wireless communications, we can observe quantifiable
system outputs evolving in time (an essential property of a dynamical
system) at every level of the ISO model hierarchy, where some aspects
of the behavior are seen as governed by simple rules and very ordered,
while others are seemingly random. The random cases attract more attention
as they are more difficult to predict and harness, and they usually
include more interaction with the external factors: users and the
environment. The question of whether to model these parts as dynamical
systems is a multi-layered one.}

\textcolor{black}{Not all systems can be described by a low-order,
low-complexity model. We could model the universe using elementary
particles positions and momenta at a certain point in time, but the
computational power for such a variant of Laplace's demon surpasses
the information limits of the same universe greatly. Measurements
might be hard to obtain as well, and while throwing a die is a perfectly
deterministic process with simple dynamics, not knowing the precise
(and ever-changing) initial conditions at the time of the throw makes
it look random }\textcolor{black}{(the level of understanding of the
processes and the granularity of the model ties closely with the traditional
understanding of uncertainty in engineering models \cite{kabir2018neural,kabir2021uncertainty})}\textcolor{black}{.
While randomness is usually portrayed as the opposite of determinism
and modeled in a different fashion, dynamical systems embrace it as
much as they embrace determinism. Random dynamical systems are a generalized
version of the deterministic ones, as they allow for a stochastic
component. In wireless communications, this is the channel equation
we start with: a linear random dynamical system with noise as a random
component.}

\textcolor{black}{The history of dynamical systems theory shows its
early bond with wireless communications, as some of the fundamental
dynamical systems theory concepts introduced by Poincare came from
his wireless telegraphy seminars (1908). The development of both disciplines
in the following decades continued going hand in hand with the oscillators
that were the central object of interest in dynamical systems theory
and an invaluable component of every radio device from their inception.
The perspectives we investigate here are those opened in dynamical
systems theory once the theory of deterministic chaos was established
in the second half of the 20th century, with the notions of sensitive
dependence on the initial conditions, fractal dimension attractors,
ergodicity, etc. The path we chose is one of understanding the already
existing dynamical phenomena within wireless systems and putting them
to use. }

\textcolor{black}{The other road in observing dynamical systems in
wireless communications is the one of chaotic communications. This
has been the dominant topic in the area since the early nineties when
the possibility of synchronization of two chaotic systems was demonstrated
\cite{pecora1990synchronization}.  This topic has repeatedly been
surveyed in the past and it represents the chaos added to a communication
system, not the one found existing within it.}

\textcolor{black}{We begin our story by presenting the basics of dynamical
systems theory. This will help us appreciate the efforts made in the
past to identify elements of dynamical systems in wireless communications
settings. These efforts will then be presented systematically, mapped
onto the media layers of ISO OSI }\textcolor{black}{(International
Organization for Standardization - Open Systems Interconnection)}\textcolor{black}{{}
model to put the concepts into a context and to suggest the ways to
proceed with the research today. Then we proceed with our central
claim: the combination of dynamical systems theory and machine learning
has a potential to radically change wireless communications performance.
We offer some initial results motivated by recent developments in
the field as a motivation for further work.}

\section{\textcolor{black}{Dynamical systems and chaos}}

\textcolor{black}{While chaotic behavior remains the trademark of
dynamical systems theory and the most interesting exhibit in its zoo,
it also remains a rare catch. When we aim to understand a dynamical
system, we are interested in its stability, periodicity, controllability,
observability: properties of the system acting on its own and under
our influence.}

\textcolor{black}{Giving attributes of a dynamical system to signal
components previously considered to be random noise: (1) allows a
better prediction, which in turn enables better suppression of interference;
(2) opens an opportunity to examine it as ``a feature, not a bug''-
i.e. adds another degree of freedom; (3) offers a physical interpretation.}

\subsection{\textcolor{black}{Signals, Phase Space, Attractors}}

\textcolor{black}{A dynamical system is, once we know all of its degrees
of freedom and sources of dynamics, a system of differential or difference
equations depending on whether we work in continuous or discrete time.
However, we tend to know so much only about very simple models seen
in nature, or the models we devise ourselves. Usually, a dynamical
system seen in the wild is a black box for us.}

\textcolor{black}{Both the system of equations and a black box take
inputs, change their states and produce outputs, which all change
in time. The number of state variables of a system is its }\textcolor{black}{\emph{order}}\textcolor{black}{,
the order of the equations' system in case we have a mathematical
description. The outputs are usually some of the states of the system
visible to us. }

\textcolor{black}{Traditionally, system identification, i.e. building
a model from the limited knowledge about the black box, is a matter
of statistics and special sets of test inputs. This is the way the
wireless channel is estimated with a (pilot) signal. Often, we try
to obtain static or linear dynamic models as they are easy to work
with. However, they are usually valid only within a narrow time or
parameter interval. Nonlinear systems ask for different identification
and modeling methods.}

\begin{figure}
\begin{centering}
\textcolor{black}{\includegraphics[width=1\columnwidth]{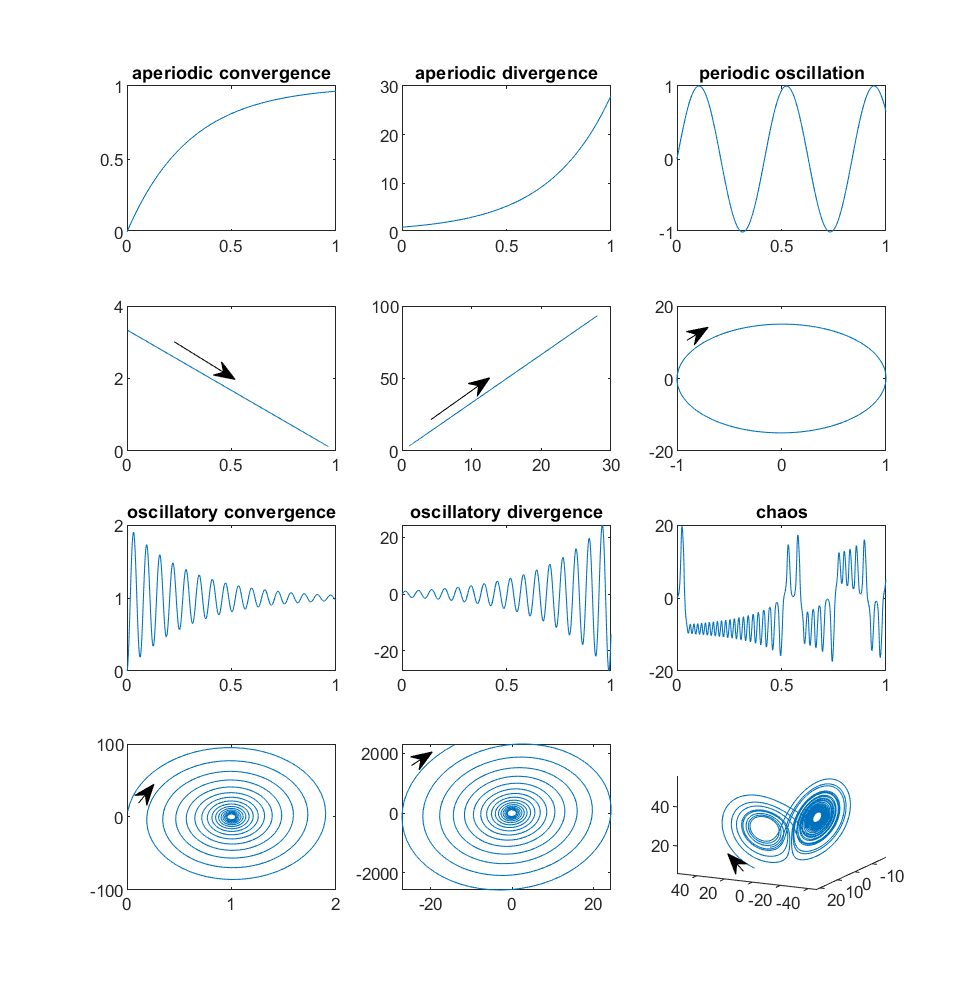}}
\par\end{centering}
\textcolor{black}{\caption{Signals in time (odd rows) and corresponding attractors in phase space
(even rows). \textcolor{black}{The arrow indicates the direction along
which system evolves in phase space (from start to end)}.}
\label{attractors}}

\end{figure}

\textcolor{black}{Having $n$ state variables, it is often useful
to plot them in an $n$-dimensional coordinate system, the }\textcolor{black}{\emph{phase
space}}\textcolor{black}{. This is where phase trajectories are observed;
typical examples are shown in Figure \ref{attractors}. The unstable
systems may diverge to infinity either quasi-periodically or aperiodically;
the stable systems may converge to an equilibrium in the same manner;
while the periodical systems remain confined to a cycle. However,
the chaotic systems were found not to follow any of these patterns:
they end up confined in a bounded part of state space called the }\textcolor{black}{\emph{attractor}}\textcolor{black}{{}
(unlike unstable systems), traverse it in a non-periodic manner (unlike
periodic systems) and never converge to a single equilibrium (unlike
stable systems). Often it is an example of a motion around two equilibria
and jumping from orbiting one to orbiting the other, as seen in the
celebrated Lorenz attractor shown in Fig. \ref{attractors}. The evolution
of trajectories on the attractor demonstrates two important properties
of a chaotic system. The first one is the }\textcolor{black}{\emph{sensitive
dependence on initial conditions}}\textcolor{black}{{} as two arbitrarily
close phase space trajectories will separate exponentially fast on
the attractor, rendering prediction of future motion impossible in
the long run. The second one is }\textcolor{black}{\emph{ergodicity}}\textcolor{black}{:
a phase trajectory will get arbitrarily close to any point on the
attractor, given enough time has passed. }

\begin{figure}
\begin{centering}
\textcolor{black}{\includegraphics[width=1\columnwidth]{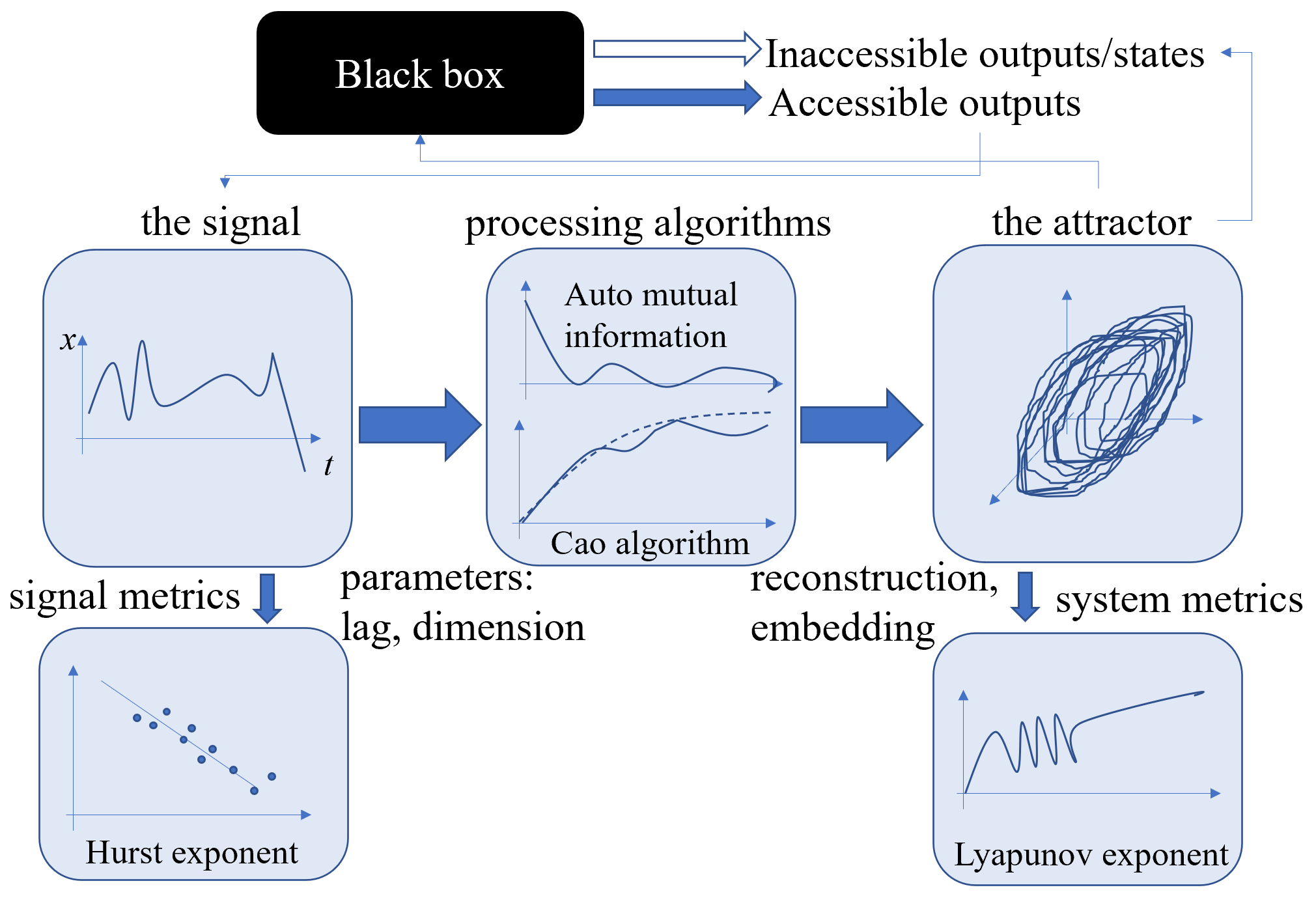}}
\par\end{centering}
\textcolor{black}{\caption{\textcolor{black}{The tool-chain for dynamical systems and chaos.
The attractor reconstruction is based on the Takens theorem, which
shows that it is possible to reconstruct an attractor based on a single
output signal and delayed versions of it serving as the other state
variables. The delay (lag) that should be used for the remaining state
variables can be determined based on the auto-mutual information function
of the signal. The use of auto-mutual information function is another
hint of how closely intertwined dynamical systems theory and information
theory are. After determining the delay and applying an algorithm
to determine the dimension the attractor is going to be embedded into,
i.e. the order of the system, an attractor can be generated. While
not perfectly the same as the original, it can reveal a lot about
system's dynamics and serve as a foundation for the calculations of
relevant metrics. The importance of this process is the transition
from signals to systems, from the particular inputs and outputs to
the general mechanism.}}
\label{tools}}

\end{figure}

\textcolor{black}{The signals and the attractors are easily obtained
when the mathematical model of the system exists: even though the
nonlinear differential/difference equations governing the dynamics
are usually not solvable in closed form, a numerical solution can
be found. However, while the system is a black box, we usually have
only few (typically, only one) system outputs available. How to reconstruct
the other state variables? How many to reconstruct in the first place?
The signal analysis and processing for chaotic dynamical systems has
a toolbox for this task. While a detailed description goes out of
the scope of this paper, Fig. \ref{tools} gives an overview of the
attractor reconstruction and quantitative analysis of the results. }

\subsection{\textcolor{black}{The Metrics}}

\textcolor{black}{Dynamical systems are quantitative, so many metrics
are devised to assess and categorize them. From the viewpoint of control
and stability, measures of stability margins describe how stable and
robust a system is and how much disturbance it could take without
a failure. }

\textcolor{black}{For chaotic systems, the measure called the }\textcolor{black}{\emph{Lyapunov
exponent}}\textcolor{black}{{} gained importance. In a conventional,
non-chaotic deterministic system, it is a negative exponent which
describes how fast two separate phase trajectories converge. In chaotic
systems, however, we have already learned that even infinitesimally
close trajectories diverge at an exponential rate: a positive Lyapunov
exponent describes this dynamic. Positive Lyapunov exponent thus became
a symbol of chaos and a basis for its measure.}

\textcolor{black}{An $n$-dimensional dynamical system has $n$ Lyapunov
exponents, and if it is chaotic, at least one of them is positive,
resulting in sensitive dependence on the initial conditions. If at
least two Lyapunov exponents are positive, we speak of hyper-chaotic
systems. Now, how do we measure just how chaotic a system is, and
how do we compare two chaotic systems in any way?}

\textcolor{black}{One possible entropy-based approach, the }\textcolor{black}{\emph{Kolmogorov-Sinai
entropy}}\textcolor{black}{{} is related to the positive Lyapunov exponents
and may give us a hint of just how unpredictable a chaotic system
is, with a numerical value between 0 (non-chaotic deterministic systems)
and infinity (purely random systems). The link with Shannon's entropy
and information theory in general is straightforward, but interesting
non-trivial results linking Lyapunov exponents of random dynamical
systems and entropy in wireless channels suggest there is more to
it than it meets the eye \cite{holliday_capacity_2006}.}

\textcolor{black}{Another metric is the dimension. In the usual sense,
we perceive dimension of a geometric construct as an integer, living
in a 3D space, observing images as 2D projections, etc. However, if
we make a finer measure of the dimension (e.g. }\textcolor{black}{\emph{Hausdorff
dimension}}\textcolor{black}{), to describe just how much of space
the object whose dimension we measure takes, we discover that not
everything is integer-dimensional. In particular, observing the Lorenz
attractor in Figure \ref{attractors} may lead to a conclusion that
it is not exactly 2-dimensional: its dimension is in fact just slightly
larger than two. The non-integer dimension of strange attractors is
thus another indicator and a measure of chaos.}

\textcolor{black}{When speaking of non-integer dimensions, it is necessary
to mention self-similarity and fractals, as one of the most often
mentioned features in dynamical systems and chaos, at least in the
popular view. Fractals are self-similar structures, as a zoomed in
part of it looks just like the bigger one, iterating the structure.
The infinitely rough structure of a fractal results in its non-integer
Hausdorff dimension, while retaining integer topological dimension
(e.g. fractals with the topological dimension of 2 can fit in a plane).}

\section{\textcolor{black}{Up and Down the Media Layers}}

\begin{table*}
\textcolor{black}{\caption{\textcolor{black}{Examples of dynamical systems research from literature,
reflecting continuous, but not intensive efforts in the area.}}
\label{tbl}}
\begin{centering}
\textcolor{black}{\small{}}%
\begin{tabular}{|>{\centering}p{0.3\textwidth}|>{\centering}p{0.6\textwidth}|}
\hline 
\textcolor{black}{\small{}Reference} & \textcolor{black}{\small{}Description}\tabularnewline
\hline 
\textcolor{black}{\small{}Holliday }\textcolor{black}{\emph{\small{}et
al. }}\textcolor{black}{\small{}2006, \cite{holliday_capacity_2006}} & \textcolor{black}{\small{}Tackling the open problem of capacity analysis
of channels characterized as Markov chains by interpreting it in terms
of Lyapunov exponents}\tabularnewline
\hline 
\textcolor{black}{\small{}Tannous }\textcolor{black}{\emph{\small{}et
al}}\textcolor{black}{\small{}. 1991, \cite{tannous_strange_1991}} & \textcolor{black}{\small{}Low-order chaos in the multipath propagation
channel with a strange attractor having a dimension between 4 and
5}\tabularnewline
\hline 
\textcolor{black}{\small{}Galdi }\textcolor{black}{\emph{\small{}et
al. }}\textcolor{black}{\small{}2005, \cite{galdi_wave_2005}} & \textcolor{black}{\small{}Accessible survey of ray chaos}\tabularnewline
\hline 
\textcolor{black}{\small{}Sirkeci-Mergen \& Scaglione 2005, \cite{sirkeci-mersen_continuum_2005}} & \textcolor{black}{\small{}Continuum approximation of dense wireless
networks, making a continuous distribution of nodes}\tabularnewline
\hline 
\textcolor{black}{\small{}Costamagna }\textcolor{black}{\emph{\small{}et
al.}}\textcolor{black}{\small{} 1994, \cite{costamagna_channel_1994}} & \textcolor{black}{\small{}First in a decade-spanning series of works
suggesting a possible way of having an attractor from a known chaotic
mapping to represent a channel error model. As the next order of approximation,
one may use several attractors stitched together to get a behavior
closer to the dynamics observed in experiments, and different channels
might be represented with different system parameters. }\tabularnewline
\hline 
\textcolor{black}{\small{}Savkin }\textcolor{black}{\emph{\small{}et
al.}}\textcolor{black}{\small{} 2005, \cite{savkin_medium_2005}} & \textcolor{black}{\small{}MAC for wireless networks as a hybrid system.
The discrete flow of data (fast dynamics) in the network is approximated
as a continuous process, while the status of nodes in the network
(working/not working, slow dynamics) remains a discrete time variable.}\tabularnewline
\hline 
\textcolor{black}{\small{}Scutari}\textcolor{black}{\emph{\small{}
et al}}\textcolor{black}{\small{}. 2008, \cite{scutari_competitive_2008}} & \textcolor{black}{\small{}Applying tools from dynamical systems theory,
stability theory and nonlinear control for multiuser MIMO (multiple
input multiple output) systems}\tabularnewline
\hline 
{\small{}Moioli }\emph{\small{}et al. }{\small{}2020, \cite{moioli2020neurosciences}} & {\small{}Drawing inspiration from dynamical system inside human brain
for sixth generation communication systems, and seeking ways to integrate
those communication systems with the human brain}\tabularnewline
\hline 
\end{tabular}{\small\par}
\par\end{centering}
\end{table*}

\textcolor{black}{The quest for chaos and other interesting dynamical
system properties was a hot topic with all features of a bandwagon
at the end of the last century. The development of algorithms described
in the previous section to determine chaoticity of the data (and the
systems generating it) brought a series of investigations and results
in different areas of science and engineering. Chaos was looked for
in the phenomena previously considered random (from the stock market
to temperature oscillations), and dynamical system formulation was
sought in the fields where the systems were not considered quantitative
at all, such as learning processes and interpersonal interaction. }

\begin{figure}
\begin{centering}
\includegraphics[width=0.9\columnwidth]{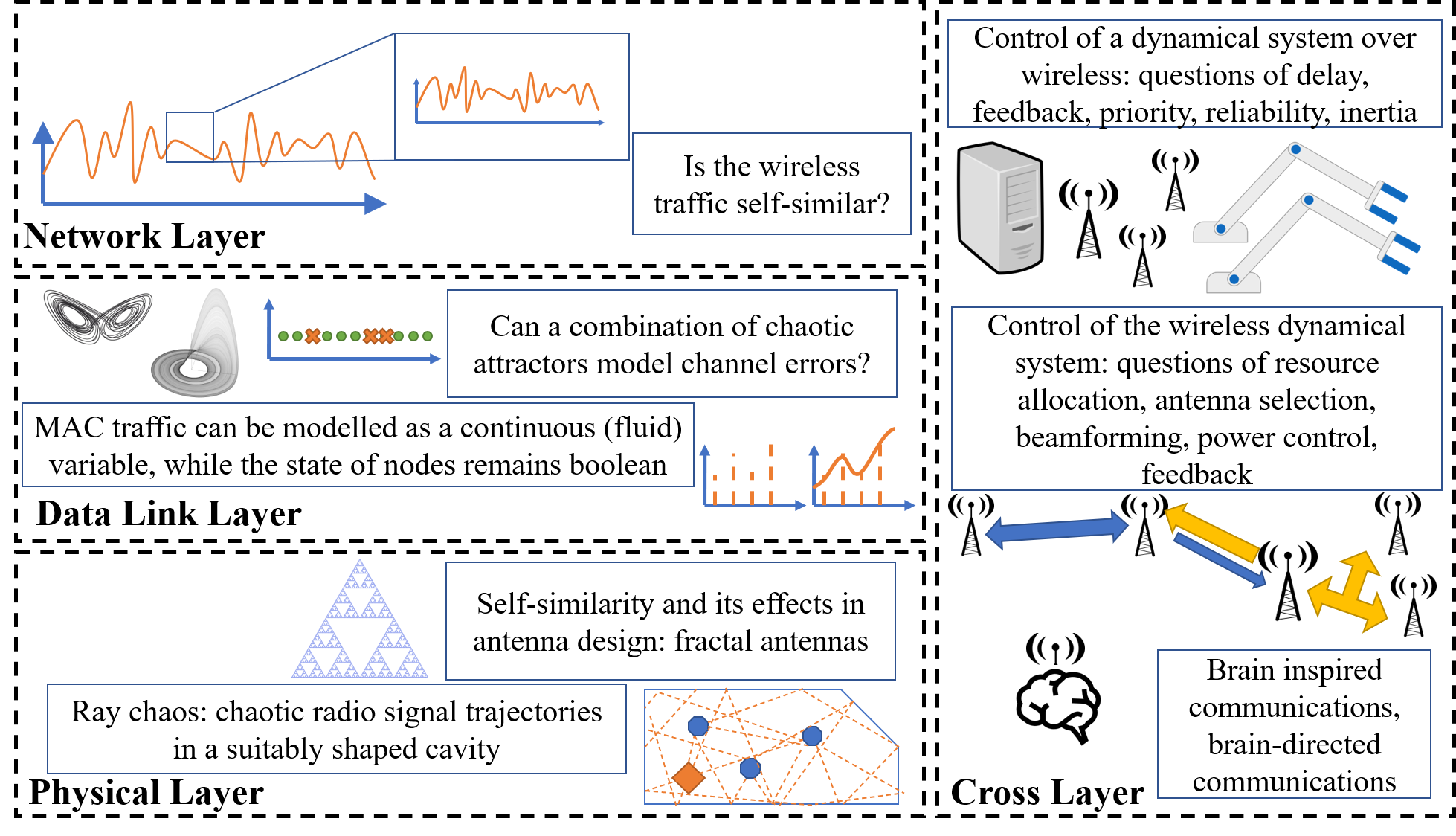}
\par\end{centering}
\textcolor{black}{\caption{What has been done in dynamical systems research over the years? Fractals,
chaos, and dynamical systems control have been put to use on different
layers of the media part of ISO OSI model.}
\label{iso-osi}}

\end{figure}

\textcolor{black}{Wireless communications followed suit, and the aspects
examined were all over the media layers of the ISO OSI model describing
it. The logic behind is simple: where there is a time series or a
signal to be measured, there is a dynamical system producing it. And
as we have seen, we can find out the details about the system from
the outputs, even more so if it is chaotic. A preview of this historical
review is given in Figure \ref{iso-osi} and elaborated upon in this
section. Illustrative references related to different layers are presented
in Table \ref{tbl}.}

\subsection{\textcolor{black}{The Physics}}

\textcolor{black}{The wireless channel has been viewed as a dynamical
system from the early beginning, but keeping a lot of its effect on
the signal under the umbrella of random noise and unpredictable changes.
Chaos-theoretic tools and new trends in dynamical systems theory have
provided some hope to distinguish the magic of noise from the science
of complexity and chaos in the channel. }

\textcolor{black}{Detecting low order chaos in multipath propagation
in the early nineties \cite{tannous_strange_1991} opened some promising
paths for the future, but similar experiments performed later show
that the deterministic chaos is not always present. The focus turned
to sea clutter in radar context, but again the results have been inconclusive,
and again coming short of the revolutionary applications chaotic sea
clutter nature might have in detection and tracking of targets near
the sea surface, if and when confirmed \cite{haykin_uncovering_2002}.}

\textcolor{black}{For a moment, examine the billiard ball system in
which a tiny difference in the initial condition drastically changes
the trajectory of the billiard ball, depending on the table shape
and the position of other balls. The billiard ball analogy works to
an extent in the context of electromagnetic waves in a chamber and
explains }\textcolor{black}{\emph{ray chaos}}\textcolor{black}{: the
divergence of two electromagnetic waves originating from almost the
same place, due to multiple bounces off the environment \cite{galdi_wave_2005}.
Ray chaos is achievable in other wave settings as well, e.g. the acoustic
case. The case of sound waves is actually the one from which we borrow
a useful application of ray chaos for radio: enabling good time reversal.}

\textcolor{black}{Time reversal in acoustics is the idea of creating
a sound wave which is a time-reversed version of a given sound wave
received: like an echo, but converging into the point of the original
wave's source. This idea \cite{fink_time_1992} extends to optics
and radio communications, serving as the core principle of several
localization and communication schemes, e.g. conjugate beamforming.
In a rich scattering environment with a lot of reflections, the bounces
in the multipath quickly map the entire space and allow good time
reversal. A chaotic environment has this property as well, as it exhibits
ergodicity: the system will eventually pass through any part of its
attractor, if we can afford to wait.}

\textcolor{black}{Time in dynamical system can either be a discrete
or a continuous variable. An interesting example of blurring the lines
between the two is the calculation of a continuous-time limit of a
discrete system. The continuous model is a limit of a discrete model
as the discrete time interval shrinks. If the systems are spatial
instead of temporal (i.e. a distance plays the role of time), the
shrinking distance produces a continuum in the limit process. In a
network, the infinitely increasing density aims for the limiting continuous
process \cite{sirkeci-mersen_continuum_2005} (Table \ref{tbl}).
The resulting continuous system may either be an approximation good
enough, or a way to obtain lower/upper bounds on the discrete system
performance.}

\subsection{\textcolor{black}{Data Link Layer}}

\textcolor{black}{Taking a step up on the ladder, we stay in the realm
of channel modeling. The question of channel errors, their occurrence,
modeling and distribution is observed and answered at the data link
layer \cite{costamagna_channel_1994}. Channel error models are dynamical
systems: they have a time evolution, have outputs in the form of discrete
time series, but are they chaotic? Can a model based on chaotic systems
encompass the dynamics of channel errors?}

\textcolor{black}{The subtle problem of modeling is bringing in all
the existing elements of the real system into the model while keeping
out all the non-existing ones. A certain dynamical system model that
happens to be chaotic may represent certain statistical properties
of the real system well. However, if the original system is not chaotic,
the model is bringing a lot of rich, but unwanted properties of its
own.}

\textcolor{black}{The example of channel errors is an example of a
discrete-time system, but the processes we describe at this layer
could be inherently a combination of both continuous and discrete
time modes. This is the case of a discrete system interacting with
a continuous one, one being embedded within the other. Medium access
control (MAC) for wireless networks can be modeled as such a system
\cite{savkin_medium_2005} (Table \ref{tbl}).}

\subsection{\textcolor{black}{Self-similar traffic across the layers }}

\textcolor{black}{The characteristics of the traffic in (wireless)
networks have been investigated for years from the perspective of
signal processing, statistics, linear and nonlinear signal theory.
The identification of statistical properties and nonlinear model parameters
of both wireless LAN }\textcolor{black}{(local area network)}\textcolor{black}{{}
and IP }\textcolor{black}{(internet protocol)}\textcolor{black}{{} traffic
would enable better models and therefore better control, prediction,
caching and routing. An elementary question is that of memory: is
the time series of the wireless LAN traffic showing any long range
dependence over time, or it is just coin tossing at every time instant?
The same question was raised at the network layer for the IP traffic,
observing the TCP }\textcolor{black}{(transmission control protocol)}\textcolor{black}{{}
congestion control. The rise and the fall of the Nile, a celebrated
example of complex patterns in nature is an example of a system's
memory. In the study of the Nile, Hurst introduced the notion of Hurst
Exponent as a measure of long term dependence in data. It was not
too surprising when the Hurst exponents were found to be directly
related with fractal dimension of self-similar data, as the memory
of the system described by the Hurst exponent is affecting the smoothness
of a signal and the fractal scale. }

\textcolor{black}{The wireless LAN traffic was expected to be self-similar:
in the mid-nineties, the self-similarity in Ethernet traffic was detected
and the research community had a field day in wireless traffic. The
packets, the bursts, the round trip time, the errors, self-similarity
appeared to be everywhere. A contributing factor for this illusion
was  the non-existence of the one right way to determine and quantify
self-similarity in data. Another reason was the existence of external
effects creating an illusion of self similarity, such as periodic
interference. The case of IP traffic was more fruitful in self-similarity
terms, but less relevant to the wireless context.}

\textcolor{black}{Down on the physical layer, some researchers have
suggested the introduction of generated self-similar signals into
wireless communication. One example is the use of self-similar carriers
(and consequently, fractal modulation), which had limited adoption
\cite{wornell_wavelet-based_1992}. The physical fractal structure,
however, has a long tradition in antenna design for wireless communications.
A typical fractal antenna has a fractal shape and puts to use the
two advantages of fractals, the existence of scaled structures and
the space-filling properties. The self-similar scaled structures offer
different scales of length for the antenna to work at, directly resulting
in similar effects for different wavelengths. This does not necessarily
mean good performance for the said wavelengths, so the story of fractal
antennas is not that straightforward. The space-filling property is
related to the capability of all fractals to achieve dense packing
in some parts (in terms of antennas, it is dense packing of wires
within a small surface area). }\textcolor{black}{For chaotic signals,
it was suggested that their transmission would be robust to multipath
effects and severe noise \cite{ren2013wireless}.}

\subsection{\textcolor{black}{The Cameos: brain, control and games}}

\textcolor{black}{Wireless communications have a long term relationship
with several disciplines heavily relying on dynamical systems. This
part of the story must not be overlooked, so we examine the control
over wireless and the game theory in the wireless setting.}

\textcolor{black}{First, we note the cameo role of dynamical systems
in wireless communication as seen in differential games. Defined as
a game over a dynamical system, a differential game can model, control
and optimize various processes in wireless communications. Again,
it is very often a hybrid dynamical system the control is performed
on, dealing with external dynamics such as drone or robotic movement,
but also the intrinsic ``dynamics'' of wireless including power
control, transmission rates and delays \cite{scutari_competitive_2008}.
To make it realistic and useful for practical considerations, the
models of the dynamical systems have to be faithful to the reality.
This means that the differential game theory eagerly awaits the results
of everything dynamical systems research can get from the wireless
of the day. }

\textcolor{black}{Closely related to game theory is the control over
wireless: it uses the wireless channel as the control signal medium
and has to deal with all of its peculiarities, }\textcolor{black}{determinism
and randomness alike (a timely example is that of teleoperation \cite{kebria2018control}}\textcolor{black}{).
The control engineers got rid of the wires cluttering the factory
and decreasing mobility, but had to face a whole new world of wireless
communications, technologies and protocols. The distributed control
system just got another dynamical system on top of it, between its
nodes: the wireless. The ubiquitous wireless sensor networks are essentially
control networks, just without actuators. And there again the dynamical
systems emerge: the ones whose outputs are measured by the sensors
and the ones the sensory data travels through. }

Finally, one area of wireless communications that recognises the importance
of dynamical systems and chaos is at the crossroads of engineering
and neuroscience, where researchers both draw inspiration from the
dynamical systems and chaotic activity in brain for bio-inspired network
and device design, and seek ways to integrate modern wireless communications
with human organism (brain in particular), for health-related applications
of the next generation communications \cite{moioli2020neurosciences}.

\section{Learning Dynamical Systems in Wireless Communications}

\begin{figure}
\begin{centering}
\includegraphics[width=6.5cm]{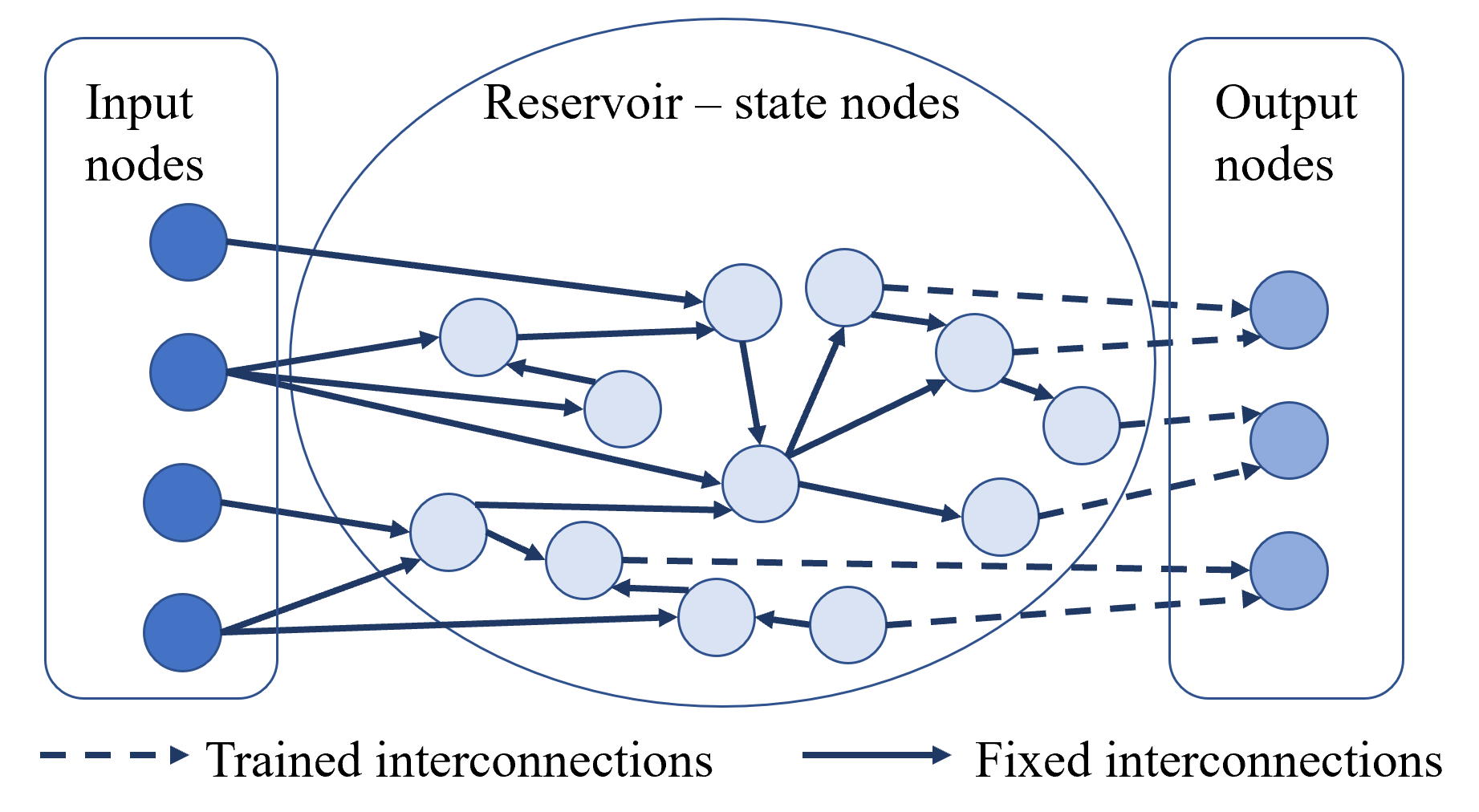}
\par\end{centering}
\caption{Reservoir computing principle. Borrowing the structure from neural
networks, reservoirs are dynamical systems with temporal dynamics
and random fixed connections between state nodes (randomly selected
when the reservoir is created and fixed for the future, unchanged
during training). Only the output connections (toward the output nodes)
are the result of training, which reduces the complexity of learning
and allows the reservoir computer to use the large dynamic reservoir
for different dynamical systems.}
\label{rezer}

\end{figure}

\textcolor{black}{The two decades of the new century saw the new artificial
intelligence spring, growing data availability and the growing capacity
for data handling. In wireless communications, an early major milestone
was the pioneering work on reservoir computing \cite{jaeger_harnessing_2004}.
This concept, illustrated in Fig. \ref{rezer} is the quintessential
machine learning model for dynamical systems: the intermediate stage
between inputs and outputs in this network is a dynamical system on
its own\textendash it does not attempt to adapt to the actual system,
except for the subset of connections leading to the outputs. As such
a general dynamical systems tool, reservoir computing quickly left
the realm of wireless communication }\textcolor{black}{(only two reservoir
computing applications cited in the most recent surveys of learning
in this area \cite{luong2019applications,jagannath2019machine})}\textcolor{black}{{}
and found greener pastures. It aims, together with other machine learning
paradigms, at providing a helping hand in predicting the behavior
of otherwise hard-to-anticipate nonlinear systems, but we are still
arguably waiting for revolutionary results. They are within the reach
once we offer a helping hand to the machine learning as well: we need
better models of dynamical systems for it to work on. More detailed
models may ask for more computing power, but we do have it now, and
in turn they greatly reduce the search space for the machine learning
and allow it to focus its efforts.}

\begin{figure}
\begin{centering}
\subfloat[\textcolor{black}{The sum rates increase over 20\%, latency drops
by 20\% and the computation burden (directly related to the power
consumption reduction and minimization of hardware) drops by 90\%
in scenarios observed by increasing the prediction horizon from one
unit interval at which prediction is currently possible (one Lyapunov
time) to eight unit intervals. The conversion of prediction time gain
into gain in other metrics was performed in straightforward manner.}]{\begin{centering}
\textcolor{black}{\includegraphics[width=0.8\columnwidth]{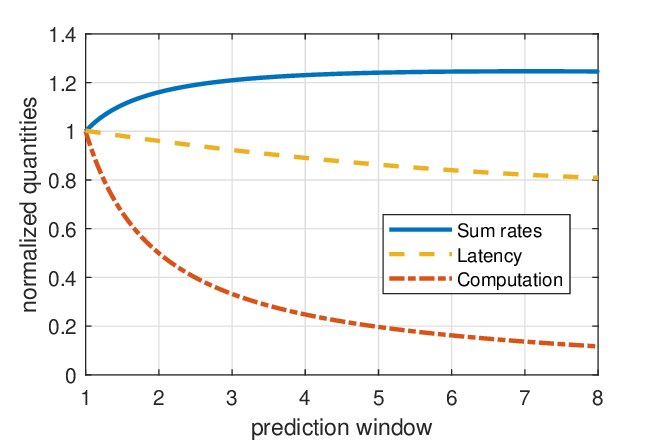}}
\par\end{centering}

}
\par\end{centering}
\subfloat[\textcolor{black}{The sum rates in a fast changing system with many
users in a cell depend on up-to-date channel state information (CSI)
for a large number of users (our scenario is a co-located Massive
MIMO base station with users in motion). Modelling the effect of CSI
acquisition at the beginning of each coherence interval on the sum
rates and the gradual detoriation of the CSI relevance over increasingly
larger intervals leads to the expected results of sum rate increase
by the virtue of less frequent need for corrections via CSI acquisition.}]{\begin{centering}
\textcolor{black}{\includegraphics[width=0.27\columnwidth]{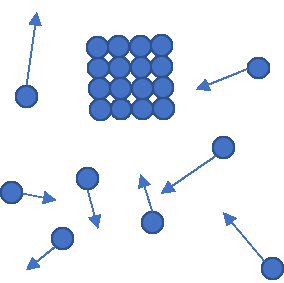}}
\par\end{centering}

}\textcolor{black}{\hfill{}}\subfloat[\textcolor{black}{Distributed massive MIMO with antenna selection
also relies heavily on constant CSI updates for a large number of
users. Using the power of prediction increase from 1 to 8 coherence
intervals and at the same time allowing the vast network of transmitters
and receivers to take a continuous form in its limit, we trade the
burden of CSI acquisition and consequent optimisation for the burden
of dynamical system prediction, which reduced the total computation.}]{\begin{centering}
\textcolor{black}{\includegraphics[width=0.27\columnwidth]{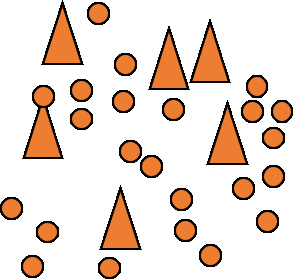}}
\par\end{centering}
}\textcolor{black}{\hfill{}}\subfloat[\textcolor{black}{We observe the time needed to deliver a message
in a multi-hop system (somewhat similar to the dynamics of epidemic
spread, a famous dynamical system) as a measure of latency. The message
is transferred from the initial point with the goal of reaching the
fixed destination via relays which move on a walk assumed predictable
by the reservoir computer. Before handing the message from one relay
to another, it is checked whether the new carrier will come closer
to the destination than the original carrier within the prediction
horizon.}]{\begin{centering}
\textcolor{black}{\includegraphics[width=0.27\columnwidth]{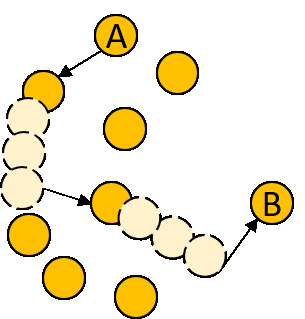}}
\par\end{centering}
}

\textcolor{black}{\caption{\textcolor{black}{The direct improvements of communication quality
based on better dynamical system prediction: (a) sum rate, computational
burden and latency improvements, (b) co-located massive MIMO scenario
for sum rate analysis, (c) distributed massive MIMO scenario for computation
analysis, (d) multi-hop scenario for latency analysis.}}
\label{rezovi}}
\vspace{-0.5cm}
\end{figure}

\textcolor{black}{Recent developments suggest that even the black
box approach in which the model has no knowledge of the actual physicality
of the process has a lot to offer for nonlinear dynamic systems. The
application of reservoir computing to spatiotemporal chaotic systems
allowed an expansion of prediction horizon\textendash as we stated
earlier, prediction of chaotic behavior is hard as two infinitesimaly
close trajectories diverge exponentially fast, and they separate within
the time period called }\textcolor{black}{\emph{Lyapunov time }}\textcolor{black}{(cf.
Lyapunov exponents). Pathak et al. \cite{pathak_model-free_2018}
report the extension of reliable prediction window from one Lyapunov
time interval to eight Lyapunov times for a particular chaotic system
(Kuramoto-Sivashinsky equation). This promising result, applicable
to wider classes of chaotic systems, motivates our investigation of
its effects in wireless communications: how to convert the information
about the future into gains in basic communication quality parameters
such as sum rates and latency? Improvement in these metrics has been
a significant driver of the technological progress, a major argument
for inclusion of new approaches (massive MIMO, mmWave) in the new
standards, and a defining aspect of the current wireless communications
generation, the 5G. In this analysis, we focus on }\textcolor{black}{\emph{sum
rate increase}}\textcolor{black}{, }\textcolor{black}{\emph{computational
burden decrease }}\textcolor{black}{and }\textcolor{black}{\emph{latency
decrease}}\textcolor{black}{{} (Fig. \ref{rezovi}(a)). These results
from simple use case scenarios (Fig. \ref{rezovi}(b-d)) suggest significant
direct benefit from dynamical systems approach to system state prediction
and add quantitative incentives to our initiative for dynamical systems
research: even relatively short prediction has a major effect on the
observed variable, justifying investing resources into it.}

The scenarios are inspired by use cases on next generation networks.
The co-located massive MIMO example (\textcolor{black}{\ref{rezovi}(b))
explores the question of what benefit for sum rates can we observe
if we can predict changes in channel state information (CSI) between
coherence intervals (periods in which CSI can be considered known
and constant). }The distributed massive MIMO example (\textcolor{black}{\ref{rezovi}(c))
was used to compare the computational overhead for prediction against
the computation needed for sensing, processing, and re-calculation
of CSI updates. Finally, the mobile multi-hop system example }(\textcolor{black}{\ref{rezovi}(d))
was set up inspired by traffic scenarios (e.g. autonomous cars in
urban settings) to see how predictability of motion helps in minimising
the message delivery time. These examples and their spatiotemporal
variability link well with the original work of Pathak et al. \cite{pathak_model-free_2018}
.}

When we speak of models, their integration with machine learning,
and reservoir computing, it is interesting to note that Pecora and
Carroll, the same researchers who have founded the field of chaotic
communications with their work on synchronisation of circuits \cite{pecora1990synchronization},
have also recently done major work on understanding the effect choice
of network connections within the reservoir \cite{carroll2019network}.
Depending on the dynamical system the reservoir computer is trying
to predict the behaviour of, the structure of computer itself does
in fact matter: flipping some of the fixed connections in the reservoir
(full line arrows in Fig. \ref{rezer}) deliver different prediction
results for the same inputs. Design of structures, given information
about application domain, hence becomes an important segment of work
and asks for understanding of complex networks. Network science and
study of complex networks are additional promising techniques for
the telecommunications community's emergent toolbox of the future. The complexity is both a blessing and a curse: unintended consequences may lead to massive failures of complex networks of "smart agents" \cite{nardelli2018smart}, but carefully designed solutions have a lot of promise for the future of technology facing climate emergency and the demise of current economic systems \cite{nardelli2021virtual}.

\section{\textcolor{black}{Conclusions}}

\textcolor{black}{The time is now: the toolbox for dynamical systems
has been reinforced with machine learning techniques born out of dynamical
systems, and the dynamics of wireless communications offer much more
to work with every day. This does not mean that the work done in the
past was not important; we need to revisit it with techniques and
technology we have today, and the results might have an application
that could not have been foreseen decades ago.}

\textcolor{black}{Again, to repeat the statement we began with: the
question is not whether we can treat everything in a wireless network
as a dynamical system, but whether we can afford to do so. The demand
is high, as the 6G and the generations of wireless to follow will
benefit from getting to know the nature of the complex, dynamical
world they are creating and embedding into at the same time. It is
hard to find a use case of next generation networks where we could
not see a nonlinear differential equation waiting to be modeled: be
it the ``wetware'' integration with dynamics of a human body, motion
in the high mobility scenario, or the myriad of rays running into
potentially ray-chaotic states in Massive MIMO, terahertz communications,
metasurfaces, the dynamical systems are within reach.}



\begin{thebibliography}{1}
\providecommand{\url}[1]{#1}
\csname url@samestyle\endcsname
\providecommand{\newblock}{\relax}
\providecommand{\bibinfo}[2]{#2}
\providecommand{\BIBentrySTDinterwordspacing}{\spaceskip=0pt\relax}
\providecommand{\BIBentryALTinterwordstretchfactor}{4}
\providecommand{\BIBentryALTinterwordspacing}{\spaceskip=\fontdimen2\font plus
\BIBentryALTinterwordstretchfactor\fontdimen3\font minus
  \fontdimen4\font\relax}
\providecommand{\BIBforeignlanguage}[2]{{%
\expandafter\ifx\csname l@#1\endcsname\relax
\typeout{** WARNING: IEEEtran.bst: No hyphenation pattern has been}%
\typeout{** loaded for the language `#1'. Using the pattern for}%
\typeout{** the default language instead.}%
\else
\language=\csname l@#1\endcsname
\fi
#2}}
\providecommand{\BIBdecl}{\relax}
\BIBdecl


\bibitem{kabir2018neural}
H.~D. Kabir, A.~Khosravi, M.~A. Hosen, and S.~Nahavandi, ``Neural network-based
  uncertainty quantification: A survey of methodologies and applications,''
  \emph{IEEE access}, vol.~6, pp. 36\,218--36\,234, 2018.

\bibitem{kabir2021uncertainty}
H.~D. Kabir, A.~Khosravi, S.~K. Mondal, M.~Rahman, S.~Nahavandi, and R.~Buyya,
  ``Uncertainty-aware decisions in cloud computing: Foundations and future
  directions,'' \emph{ACM Computing Surveys (CSUR)}, vol.~54, no.~4, pp. 1--30,
  2021.

\bibitem{pecora1990synchronization}
L.~M. Pecora and T.~L. Carroll, ``Synchronization in chaotic systems,''
  \emph{Physical review letters}, vol.~64, no.~8, p. 821, 1990.

\bibitem{holliday_capacity_2006}
T.~Holliday, A.~Goldsmith, and P.~Glynn, ``Capacity of {Finite} {State}
  {Channels} {Based} on {Lyapunov} {Exponents} of {Random} {Matrices},''
  \emph{IEEE Transactions on Information Theory}, vol.~52, no.~8, pp.
  3509--3532, Aug. 2006.

\bibitem{tannous_strange_1991}
C.~Tannous, R.~Davies, and A.~Angus, ``Strange attractors in multipath
  propagation,'' \emph{IEEE Transactions on Communications}, vol.~39, no.~5,
  pp. 629--631, May 1991.

\bibitem{galdi_wave_2005}
V.~Galdi, I.~M. Pinto, and L.~B. Felsen, ``Wave propagation in ray-chaotic
  enclosures: paradigms, oddities and examples,'' \emph{IEEE Antennas and
  Propagation Magazine}, vol.~47, no.~1, pp. 62--81, Feb. 2005.

\bibitem{sirkeci-mersen_continuum_2005}
B.~Sirkeci-Mersen and A.~Scaglione, ``A continuum approach to dense wireless
  networks with cooperation,'' in \emph{Proceedings {IEEE} 24th {Annual}
  {Joint} {Conference} of the {IEEE} {Computer} and {Communications}
  {Societies}.}, vol.~4, Mar. 2005, pp. 2755--2763 vol. 4.

\bibitem{costamagna_channel_1994}
E.~Costamagna and A.~Schirru, ``Channel error models derived from chaos
  equations,'' in \emph{Proceedings of {IEEE} {International} {Conference} on
  {Systems}, {Man} and {Cybernetics}}, vol.~1, Oct. 1994, pp. 577--581 vol.1.

\bibitem{savkin_medium_2005}
A.~V. Savkin, A.~S. Matveev, and P.~B. Rapajic, ``The medium access control
  problem for wireless communication networks modelled as hybrid dynamical
  systems,'' \emph{Nonlinear Analysis: Theory, Methods \& Applications},
  vol.~62, no.~8, pp. 1384--1393, Sep. 2005.

\bibitem{scutari_competitive_2008}
G.~Scutari, D.~P. Palomar, and S.~Barbarossa, ``Competitive {Design} of
  {Multiuser} {MIMO} {Systems} {Based} on {Game} {Theory}: {A} {Unified}
  {View},'' \emph{IEEE Journal on Selected Areas in Communications}, vol.~26,
  no.~7, pp. 1089--1103, Sep. 2008.

\bibitem{moioli2020neurosciences}
R.~C. Moioli, P.~H. Nardelli, M.~T. Barros, W.~Saad, A.~Hekmatmanesh,
  P.~G{\'o}ria, A.~S. de~Sena, M.~Dzaferagic, H.~Siljak, W.~Van~Leekwijck,
  D.~Carrillo, and S.~Latre, ``Neurosciences and 6g: Lessons from and needs of
  communicative brains,'' \emph{IEEE Communications Surveys \& Tutorials},
  2021.

\bibitem{haykin_uncovering_2002}
S.~Haykin, R.~Bakker, and B.~W. Currie, ``Uncovering nonlinear dynamics-the
  case study of sea clutter,'' \emph{Proceedings of the IEEE}, vol.~90, no.~5,
  pp. 860--881, May 2002.

\bibitem{fink_time_1992}
M.~Fink, ``Time reversal of ultrasonic fields. {I}. {Basic} principles,''
  \emph{IEEE Transactions on Ultrasonics, Ferroelectrics, and Frequency
  Control}, vol.~39, no.~5, pp. 555--566, Sep. 1992.

\bibitem{wornell_wavelet-based_1992}
G.~W. Wornell and A.~V. Oppenheim, ``Wavelet-based representations for a class
  of self-similar signals with application to fractal modulation,'' \emph{IEEE
  Transactions on Information Theory}, vol.~38, no.~2, pp. 785--800, Mar. 1992.

\bibitem{ren2013wireless}
H.-P. Ren, M.~S. Baptista, and C.~Grebogi, ``Wireless communication with
  chaos,'' \emph{Physical Review Letters}, vol. 110, no.~18, p. 184101, 2013.

\bibitem{kebria2018control}
P.~M. Kebria, H.~Abdi, M.~M. Dalvand, A.~Khosravi, and S.~Nahavandi, ``Control
  methods for internet-based teleoperation systems: A review,'' \emph{IEEE
  Transactions on Human-Machine Systems}, vol.~49, no.~1, pp. 32--46, 2018.

\bibitem{jaeger_harnessing_2004}
\BIBentryALTinterwordspacing
H.~Jaeger and H.~Haas, ``\BIBforeignlanguage{en}{Harnessing {Nonlinearity}:
  {Predicting} {Chaotic} {Systems} and {Saving} {Energy} in {Wireless}
  {Communication}},'' \emph{\BIBforeignlanguage{en}{Science}}, vol. 304, no.
  5667, pp. 78--80, Apr. 2004. [Online]. Available:
  \url{http://science.sciencemag.org/content/304/5667/78}
\BIBentrySTDinterwordspacing

\bibitem{luong2019applications}
N.~C. Luong, D.~T. Hoang, S.~Gong, D.~Niyato, P.~Wang, Y.-C. Liang, and D.~I.
  Kim, ``Applications of deep reinforcement learning in communications and
  networking: A survey,'' \emph{IEEE Communications Surveys \& Tutorials},
  vol.~21, no.~4, pp. 3133--3174, 2019.

\bibitem{jagannath2019machine}
J.~Jagannath, N.~Polosky, A.~Jagannath, F.~Restuccia, and T.~Melodia, ``Machine
  learning for wireless communications in the internet of things: A
  comprehensive survey,'' \emph{Ad Hoc Networks}, vol.~93, p. 101913, 2019.

\bibitem{pathak_model-free_2018}
J.~Pathak, B.~Hunt, M.~Girvan, Z.~Lu, and E.~Ott, ``Model-free prediction of
  large spatiotemporally chaotic systems from data: a reservoir computing
  approach,'' \emph{Physical review letters}, vol. 120, no.~2, p. 024102, 2018.

\bibitem{carroll2019network}
T.~L. Carroll and L.~M. Pecora, ``Network structure effects in reservoir
  computers,'' \emph{Chaos: An Interdisciplinary Journal of Nonlinear Science},
  vol.~29, no.~8, p. 083130, 2019.

\bibitem{nardelli2018smart}
P.~H. Nardelli and F.~K{\"u}hnlenz, ``Why smart appliances may result in a
  stupid grid: Examining the layers of the sociotechnical systems,'' \emph{IEEE
  Systems, Man, and Cybernetics Magazine}, vol.~4, no.~4, pp. 21--27, 2018.

\bibitem{nardelli2021virtual}
P.~H. Nardelli, H.~M. Hussein, A.~Narayanan, and Y.~Yang, ``Virtual microgrid
  management via software-defined energy network for electricity sharing,''
  \emph{IEEE Systems, Man, and Cybernetics Magazine}, 2021.

\end{thebibliography}

\begin{IEEEbiography}[{\includegraphics[width=1in,height=1.25in,clip,keepaspectratio]{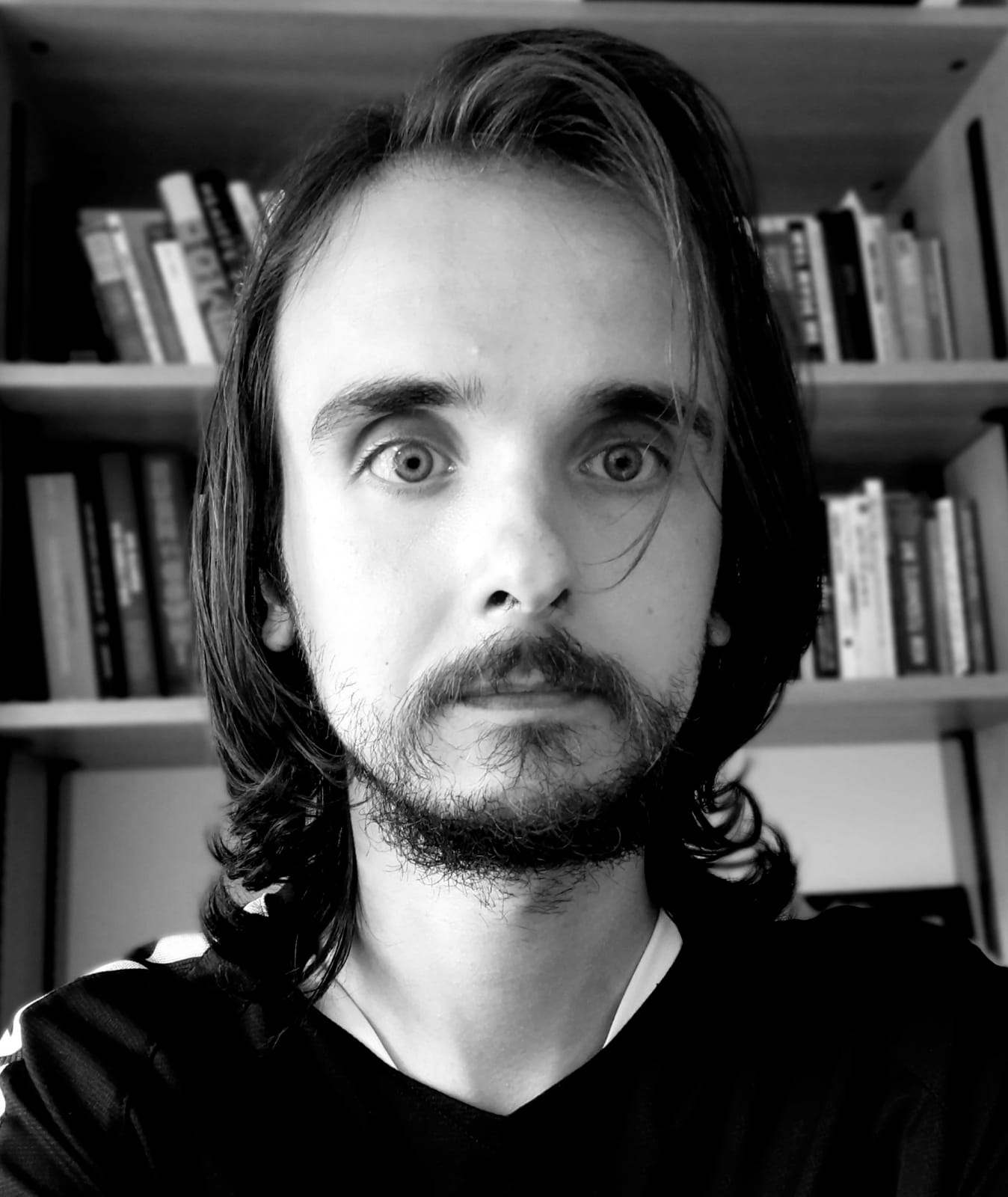}}]{Harun Siljak}

(M'15, SM'20) received the bachelor and master degrees in electronics and control engineering from the University of Sarajevo, Bosnia and Herzegovina in 2010 and 2012, respectively. He received his PhD degree in electronics and electrical engineering from the International Burch University Sarajevo, Bosnia and Herzegovina in 2015. He is currently an assistant professor in embedded systems, control and optimisation at Trinity College Dublin, Ireland, where he previously worked as a Marie Curie postdoctoral fellow (EDGE Cofund). His research interests span from complex networks and cyber-physical systems to unconventional methods of computation, communication, and control.

\end{IEEEbiography}

\begin{IEEEbiography}[{\includegraphics[width=1in,height=1.25in,clip,keepaspectratio]{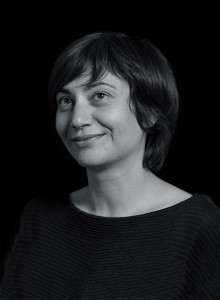}}]{Irene Macaluso}

is a Senior Research Fellow at CONNECT, Ireland's research centre for Future Networks and Communications, based at Trinity College, Dublin. Dr. Macaluso received her Ph.D. in Robotics from the University of Palermo in 2007. Dr. Macaluso's current research interests are in the area of adaptive wireless resource allocation, with particular focus on the design and analysis of market-based mechanisms in the management and operation of reconfigurable wireless networks and the application of machine learning to network resource sharing. She has published more than 80 papers in internationally peer reviewed journals and conferences and holds 2 patents. She is Executive Editor of Transactions on Emerging Telecommunication Technologies (ETT) since 2016.

\end{IEEEbiography}

\begin{IEEEbiography}[{\includegraphics[width=1in,height=1.25in,clip,keepaspectratio]{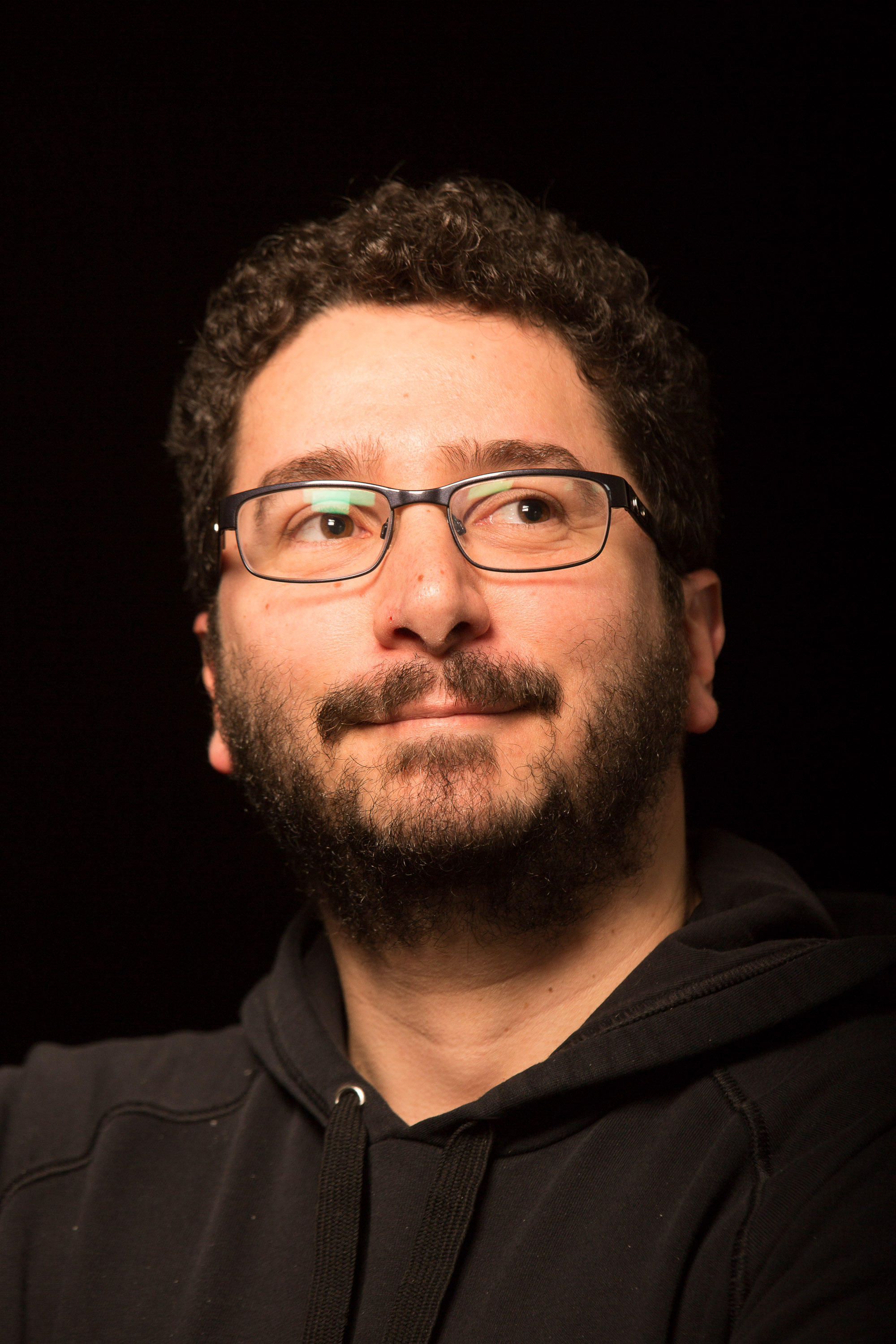}}]{Nicola Marchetti}

is Associate Professor in Wireless Communications at Trinity College Dublin, Ireland. He performs his research under the Irish Research Centre for Future Networks and Communications (CONNECT), where he leads the Wireless Engineering and Complexity Science (WhyCOM) lab. He received the PhD in Wireless Communications from Aalborg University, Denmark in 2007, and the M.Sc. in Electronic Engineering from University of Ferrara, Italy in 2003. He also holds an M.Sc. in Mathematics which he received from Aalborg University in 2010. His research interests include Self-Organizing Networks, Signal Processing for Communication, and Radio Resource Management. He has authored in excess of 150 journals and conference papers, 2 books and 8 book chapters, holds 4 patents, and received 4 best paper awards. He is a senior member of IEEE and serves as an associate editor for the IEEE Internet of Things Journal since 2018.

\end{IEEEbiography}

\pagebreak{}






\end{document}